\documentclass{amsart}
\usepackage[latin1]{inputenc}
\usepackage{amssymb,amsmath,amsthm}
\usepackage{amsfonts}
\usepackage{ifthen}
\usepackage{url}
\usepackage{paralist}

\theoremstyle{plain}
\newtheorem*{theorem}{Theorem}
\newtheorem*{reformulation}{Reformulation of the Theorem}
\newtheorem*{corollary}{Corollary}
\theoremstyle{remark}
\newtheorem{rem}{Remark}

\newcommand{\card}[1]{\ensuremath\lvert{#1}\rvert}
\begin{document}
\title[Column-partitioned matrices over rings]{Column-partitioned matrices over rings without invertible transversal submatrices}
\author{Stephan Foldes}
\address{Institute of Mathematics, Tampere University of Technology, P.O.~Box 553, FI-33101 Tampere, Finland}
\email{stephan.foldes@tut.fi}
\author{Erkko Lehtonen}
\address{Institute of Mathematics, Tampere University of Technology, P.O.~Box 553, FI-33101 Tampere, Finland}
\email{erkko.lehtonen@tut.fi}
\date{\today}
\begin{abstract}
Let the columns of a $p \times q$ matrix $M$ over any ring be partitioned into $n$ blocks, $M = [M_1, \ldots, M_n]$. If no $p \times p$ submatrix of $M$ with columns from distinct blocks $M_i$ is invertible, then there is an invertible $p \times p$ matrix $Q$ and a positive integer $m \leq p$ such that $QM = [QM_1, \ldots, QM_n]$ is in reduced echelon form and in all but at most $m-1$ blocks $QM_i$ the last $m$ entries of each column are either all zero or they include a non-zero non-unit.
\end{abstract}
\maketitle

\section{Column-partitioned matrices and transversal submatrices}

Generalizing the concept of row-reduced form of matrices over fields, we shall say that a matrix $M$ over any ring is in \emph{reduced echelon form} if among all matrices $QM$ where $Q$ is invertible it has the maximum possible number of distinct standard unit vectors appearing as columns.

\begin{theorem}
\label{mainth}
Let $R$ be any ring with identity, possibly non-commutative. Let the columns of a $p \times q$ matrix $M$ with entries in $R$ be partitioned into $n$ blocks, $M = [M_1, \ldots, M_n]$. Suppose that no $p \times p$ submatrix extracted from $M$ with columns from distinct blocks $M_i$ is invertible. Then there is an invertible $p \times p$ matrix $Q$ and a positive integer $m \leq p$ such that $QM = [QM_1, \ldots, QM_n]$ is in reduced echelon form and in all but at most $m-1$ blocks $QM_i$ the last $m$ entries of each column are either all zero or they include a non-zero non-unit.
\end{theorem}

\begin{rem}
If $n = q$ and $R$ is a field, then the Theorem is an obvious consequence of a rank-deficient matrix over a field having a null row in its reduced row-echelon form.
\end{rem}
\begin{rem}
If $n < p$, then the Theorem trivially holds with the identity matrix as $Q$ and $m = p$.
\end{rem}

In order to prove the Theorem, it will be convenient to recast it in a somewhat more general form, using the following definitions and notation for purposes of precision and simplicity in the proof.

An \emph{$A \times B$ matrix} with entries in a ring $R$ is any map $M : A \times B \to R$, where $A$, $B$ are finite sets of positive integers. The matrix product $MN$ of $M : A \times B \to R$ and $N : B \times C \to R$ is a map $A \times C \to R$ whose value on $(a,c) \in A \times C$ is defined by the usual convolution formula. For $A' \subseteq A$, $B' \subseteq B$, we denote by $M[A',B']$ the restriction of $M$ (as a map) to $A' \times B'$; thus $M = M[A,B]$. If any of $A'$ or $B'$ is a singleton $\{a\}$, then we may omit the set braces and write $a$ for $\{a\}$.

Whenever we refer to \emph{elementary row operations} on an $A \times B$ matrix $M$, we mean left multiplication of $M$ by an $A \times A$ matrix $E$ of one of the following two types:
\begin{asparaenum}
\item a diagonal matrix all whose diagonal entries are units (scaling of rows by units),
\item the sum of the identity matrix and a matrix with a single non-zero entry in an off-diagonal position (adding a multiple of a row to another row).
\end{asparaenum}
All such matrices $E$ are invertible.

For any set $B$, a \emph{partition} is a set $\Pi$ of nonempty pairwise disjoint subsets of $B$ the union of which is $B$. A \emph{partial transversal} of $\Pi$ is a subset $J$ of $B$ intersecting every partition class $K \in \Pi$ in at most one element.

\begin{reformulation}
Let $R$ be any ring with identity, possibly non-commutative, and let $M$ be an $A \times B$ matrix with entries in $R$. Consider a partition $\Pi$ of $B$ into $n$ classes, $\Pi = \{B_1, \ldots, B_n\}$. Suppose that for every partial transversal $J$ of $\Pi$ with $\card{J} = \card{A}$, the submatrix $M[A,J]$ is not invertible. Then there is an invertible $A \times A$ matrix $Q$ and a nonempty subset $A' \subseteq A$ such that $QM$ is in reduced echelon form and at most $\card{A'}-1$ of the matrices $(QM)[A',B_i]$, $1 \leq i \leq n$, can have a column containing a unit entry but no non-zero non-units.
\end{reformulation}

\section{Proof of the reformulation}

If $n < \card{A}$, then the statement clearly holds with $A' = A$. Therefore we can assume that $n \geq \card{A}$.

If there is no subset $P \subseteq B$ with $\card{P} = \card{A}$ such that $M[A,P]$ is invertible, then let $t < \card{A}$ be the largest positive integer such that there is some invertible $A \times A$ matrix $Q$ and $t$ distinct standard unit vectors that appear as columns of $QM$. (In case there is no such positive $t$, then obviously no entry of $M$ is a unit and the claimed result holds with any singleton $A'$.) Clearly $QM$ has exactly $m = \card{A} - t > 0$ rows without units, and the Theorem easily follows.

Suppose therefore that there are subsets $P \subseteq B$, $\card{P} = \card{A}$ (but none with $P$ being a partial transversal of $\Pi$) such that $M[A,P]$ is invertible. Call such subsets $P$ \emph{admissible sets.}

Define the \emph{spread} of an admissible set $P$ as the set $\{i \in \{1, \ldots, n\} : P \cap B_i \neq \emptyset\}$. The \emph{weight} (with respect to $P$) of a block $B_i$ is defined as $w_i = \card{P \cap B_i}$. The \emph{profile} of $P$ is the multiset of the weights $w_i$ where $i$ is in the spread of $P$. The \emph{profile sequence} of $P$ is the monotone increasing ordered profile of $P$. Denote the inverse of $M[A,P]$ by $Q$. For each $1 \leq i \leq n$, define the set
\[
A_i = \{r \in A : (QM)[r,c] = 1,\, c \in P \cap B_i\}.
\]

For an admissible set $P$, denote
\[
\hat{A} = \Big( \bigcup_{w_i \leq 2} A_i \Big), \qquad \hat{B} = \Big( \bigcup_{w_i \leq 1} B_i \Big),
\]
and let $D_0 = (\hat{A} \cup \{0\}) \times \hat{B} \times \hat{A}$.
The elements $(s,c,t) \in D_0$ will be considered as the \emph{arrows} of a directed graph $G_0$ with vertex set $A \cup \{0\}$, where the source of $(s,c,t)$ is $s$ and its target is $t$, and the element $c$ distinguishes between parallel arrows; we say that $(s,c,t)$ is an arrow from row $s$ to row $t$ through column $c$. For any subset $D' \subseteq D_0$, we shall mean by ``the graph $D'$'' the subgraph of $G_0$ which contains all vertices (i.e., $A \cup \{0\}$) but only those arrows that are in $D'$. Let $D$ be the set of arrows $(s,c,t) \in D_0$ that satisfy the following conditions:
$(QM)[t,c]$ is a unit;
$s \neq t$;
if $c$ belongs to a block $B_k$ of weight $1$ then $s$ is the unique member of $A_k$, else $s = 0$.
A directed path $(\alpha_1, \ldots, \alpha_l) = ((s_1, c_1, t_1), \ldots, (s_l, c_l, t_l))$ in the graph $D$ is said to be \emph{clear} if for all $1 \leq i < l$, denoting by $T_i$ the set of targets of the arrows $\alpha_j$, $j > i$, we have $QM[T_i, c_i] = 0$.

Define pairwise disjoint subsets $D_1$, $D_2$, \ldots, $D_k$, \ldots{} of $D$ inductively as follows, denoting $\bigcup_{w_i = 2} A_i$ by $T$. The members of $D_1$ are the arrows of the graph $D$ with target in $T$. The members of $D_{k+1}$ are the arrows $(s,c,t)$ of the graph $D$ whose target is the source of an arrow in the graph $D_k$ and for which there is a clear path in the graph $\bigcup_{i=1}^k D_i \cup \{(s,c,t)\}$ starting at $(s,c,t)$ and ending with an arrow with target in $T$. Let $D_P$ be the union of all $D_k$, $k \geq 1$, and call the graph $D_P$ the \emph{connection graph} of $P$, and denote it by $G_P$. We will need to distinguish two cases. If there is no directed path from $0$ to a member of $T$ in $G_P$, then we say that $P$ is of the \emph{first kind.} Otherwise we say that $P$ is of the \emph{second kind} and the length of the shortest directed path from $0$ to a member of $T$ in $G_P$ is called the \emph{connection distance} for $P$.

Let $\mathcal{P}$ be the set of all admissible sets of maximum spread (i.e., meeting as many blocks of $\Pi$ as possible). This set is quasi-ordered by the majorization relation between profile sequences. (Recall that a monotone sequence $a_1 \leq a_2 \leq \dots \leq a_s$ is said to \emph{majorize} a sequence $b_1 \leq b_2 \leq \dots \leq b_s$ when for all $1 \leq i \leq s$, $a_1 + \dots + a_i \geq b_1 + \dots + b_i$, with equality for $i = s$. Majorization is a partial order on the set of finite monotone increasing sequences of integers.) Let $\mathcal{P}_1$ be the set of maximal members of $\mathcal{P}$ (i.e., the members of $\mathcal{P}$ whose profile sequence is not strictly majorized by the profile sequence of another member of $\mathcal{P}$). Let $P$ be an admissible set in $\mathcal{P}_1$ of the first kind if such exists, else let $P$ be a member of $\mathcal{P}_1$ (necessarily of the second kind) whose connection distance is as small as possible.

\emph{Gap Condition.} We claim that there are no units in $(QM)[A_i, B_j]$ whenever $w_i \geq w_j + 2$. For, suppose, on the contrary, that $(QM)[r,c]$ is a unit for some $r \in A_i$, $c \in B_j$ with $w_i \geq w_j + 2$. Then there is a $c' \in A_i \cap P$ such that $(QM)[r,c'] = 1$, and we can make column $c$ into a standard unit vector with elementary row operations that do not affect the columns indexed by $P \setminus \{c'\}$. Thus the set $P' = P \cup \{c\} \setminus \{c'\}$ is admissible, but it either has a larger spread than $P$ (if $w_j = 0$) or it has the same spread as $P$ (if $w_j > 0$) but its profile sequence majorizes that of $P'$, a contradiction.

Since $P$ is not a partial transversal of $\Pi$, there must be a block of weight at least $2$, and there is of course a block of weight $0$. Assume first that there is no block of weight $2$. In this case, let $I = \{i : w_i > 2\}$, and the claimed result holds by choosing $A' = \bigcup_{i \in I} A_i$, because by the Gap Condition, $(QM)[A',B_i]$ does not contain a unit for any $i \notin I$, and $\card{I} < \card{A'}$.

We can thus assume that there is a block of weight $2$. If $P$ is of the first kind, let $S$ be the set of indices $i \in \{1, \ldots, n\}$ such that $w_i = 1$ and there is no arrow $(s,c,t)$ with $c \in B_i$ on any path in $G_P$ terminating in $T$. In this case we obtain the result, if we let $I = \{1, \ldots, n\} \setminus S$ and choose $A' = \bigcup_{i \in I} A_i$. For, if $i \notin I$ and $(QM)[r,c]$ is a unit for some $r \in A'$, $c \in B_i$ (such an $r$ is necessarily in an $A_k$ with $w_k = 1$: this follows from the Gap Condition for blocks of weight $0$; and if $B_i$ is a block of weight $1$ and $r \in A_k$ with $w_k = 2$, then there would be an arrow from the single element of $A_i$ to an element of $T$ through $c$, and so $i \in I$, a contradiction), then there is an $r' \in T$ such that $(QM)[r',c]$ is a non-zero non-unit.

If $P$ is of the second kind, it is clear that the connection distance is at least $2$. In $G_P$, take a shortest directed path $((s_1, c_1, t_1), \ldots, (s_l, c_l, t_l))$ from $0$ to a vertex in $T$. For the last arrow $(s_l, c_l, t_l)$ in this path, we have $t_l \in A_k$ for some $k$ with $w_k = 2$ and $c$ belongs to a block $B_j$ of weight $1$. There is a $c \in B_k \cap P$ such that $(QM)[t_l,c_l]$ is a unit and $(QM)[t_l,c] = 1$, and we can do elementary row transformations and make $B_k$ into a block of weight $1$ and $B_j$ into a block of weight $2$ with respect to a new admissible set $P' = P \cup \{c_l\} \setminus \{c\}$. These row transformations do not affect the columns indexed by $\{c_1, \ldots, c_{l-1}\} \cup P \setminus \{c\}$. Therefore $((s_1, c_1, t_1), \ldots, (s_{l-1}, c_{l-1}, t_{l-1})$ is a clear path in $G_{P'}$ from $0$ to $t_{l-1}$ and $t_{l-1}$ now belongs to the set $A_j$ of rows corresponding to a block $B_j$ of weight $2$ with respect to $P'$. The set $P'$ has the same spread, the same profile, and the same profile sequence as $P$, it is still of the second kind, but its connection distance is smaller than that of $P$, a contradiction exhausting the last possible case. This completes the proof of the Theorem.
\qed

\begin{rem}
The Gap Condition in the above proof shows that the matrix $QM$ will indeed have some rows in which some blocks are completely free of units.
\end{rem}

\section{Matrices over fields}

The Theorem above applies to any ring $R$, whether commutative or not. In the special case that $R$ is a field, the Theorem overlaps as we shall show below with Rado's matroid-theoretical generalization of Hall's theorem on systems of distinct representatives as reformulated and extended by Perfect \cite{PerfectLond,PerfectCamb}. However, the Rado--Perfect results do not apply to matrices over arbitrary rings, as the columns of such matrices do not have the combinatorial properties stipulated by matroid theory's abstract generalization of linear independence.

Perfect's version of Rado's theorem, specifically as in Theorem 2 of \cite{PerfectCamb}, states the following, when applied to any $A \times B$ matrix $M$ over a field, any partition $\Pi$ of $B$ into $n$ classes, and any positive integer $k$:
\begin{quote}
\emph{There is a partial transversal $P$ of $\Pi$ of size $k$ such that $M[A,P]$ has rank $k$ if and only if for all $\Theta \subseteq \Pi$ the rank of $M[A, \bigcup \Theta]$ is at least $k + \card{\Theta} - n$ (where $\bigcup \Theta$ denotes the union of the partition blocks in $\Theta$).}
\end{quote}
The ``only if'' part is obvious here, while the ``if'' part states in particular that if $M[A,P]$ is non-invertible for all partial transversals $P$ of $\Pi$ of size $\card{A}$ then for some $\Theta \subseteq \Pi$ the rank $\rho$ of $M[A, \bigcup \Theta]$ is less than $\card{A} + \card{\Theta} - n = \card{A} - (\card{\Pi} - \card{\Theta})$, in other words $\card{A} - \rho > \card{\Pi \setminus \Theta}$. Gaussian elimination then yields an invertible $A \times A$ matrix $Q$ and a set $A' \subseteq A$ of size $\card{A} - \rho$ such that the matrix $(QM)[A', \bigcup \Theta]$ is identically null, i.e., at most $\card{\Pi \setminus \Theta} \leq \card{A'} - 1$ of the matrices $(QM)[A', B']$, $B' \in \Pi$, can have a non-null entry.

Thus the Rado--Perfect result on independent partial transversals in matroids and our Theorem above overlap in the following Corollary, where the implication $\text{(2)} \Rightarrow \text{(1)}$ is obvious.

\begin{corollary}
Let the columns of a $p \times q$ matrix $M$ with entries in any field be partitioned into $n$ blocks, $M = [M_1, \ldots, M_n]$. The following are equivalent.
\begin{enumerate}
\item All $p \times p$ submatrices extracted from $M$ with columns from distinct blocks $M_i$ are noninvertible.
\item There is an invertible $p \times p$ matrix $Q$ and a positive integer $m \leq p$ such that in $QM = [QM_1, \ldots, QM_n]$ the last $m$ rows are null in all but at most $m-1$ blocks $QM_i$.
\end{enumerate}
\end{corollary}

The Corollary as applied to matrices with entries in the two-element field is used in \cite{Lehtonen} to establish the descending chain condition in a particular quasi-order of operations on a finite set of $k$ elements ($k \geq 3$). The quasi-ordering is based on the composition of functions from inside with the quasi-linear functions of Burle \cite{Burle}.

\section*{Acknowledgements}

The authors thank Miguel Couceiro and Jen\H{o} Szigeti for useful discussions of the subject.

\end{document}